\documentclass[11pt,a4paper]{article}
\usepackage[latin1]{inputenc}
\usepackage{amsmath}
\usepackage{amsfonts}
\usepackage{amssymb}
\usepackage{graphicx}
\usepackage{amsthm}
\newtheorem{theorem}{Theorem}
\newtheorem{lemma}{Lemma}
\theoremstyle{definition}
\newtheorem{remark}{Remark}
\newtheorem*{closeremark}{Closing Remarks}
\newtheorem*{acknolwedgment}{Author's Note and Acknowledgments}
\author{Chance Sanford}
\title{Infinite Series Involving Fibonacci Numbers Via Ap\'{e}ry-Like Formulae}
\begin{document}

\maketitle

\begin{abstract}
In this note, infinite series involving Fibonacci and Lucas numbers are derived by employing formulae similar to that which Roger Ap\'{e}ry utilized in his seminal paper proving the irrationality of $\zeta(3)$.  
\end{abstract}

\section{Introduction}

In 1979 Roger Ap\'{e}ry published his seminal paper in which he proved the irrationality of $\zeta(3) = \sum_{n=1}^{\infty}{\frac{1}{n^3}}$. In his paper, Ap\'{e}ry utilizes a technique that allows one to evaluate certain infinite series telescopically.  

In this note, we demonstrate how Ap\'{e}ry's technique may be employed to evaluate infinite series involving Fibonacci and Lucas numbers; in particular of the form $F_{n^k}$ and $L_{n^k}$ respectively.  Using Ap\'{e}ry-like formulae we are able to evaluate series that can be difficult to determine by other means.

To the best of the author's knowledge, Ap\'{e}ry's technique has not been utilized to evaluate infinite series involving Fibonacci and Lucas numbers.  In addition, the results obtained in this note are also believed to be new, except for Theorem 1, which is a new proof of an old result.\\

Before proceeding with the results, we first examine Ap\'{e}ry's formula as described in Chen \cite{chen10}.\\

For $x \neq 0$ and any sequence ${a_k}$, let

$$B_1 = \frac{1}{x}, \quad B_k = \frac{{a_1}{a_2}\dots{a_{k-1}}}{x(x+a_1)(x+a_2)\dots(x+a_{k-1})}, \quad k \geq 2$$

Then,

$$B_k - B_{k+1} = \frac{{a_1}{a_2}\dots{a_{k-1}}}{({x+a_1})({x+a_2})\dots({x+a_{k}})}$$\\

Due to telescoping, we find that:
\begin{equation*}
\frac{1}{x+a_1}+ \sum_{k=2}^{n}\frac{{a_1}{a_2}\dots{a_{k-1}}}{({x+a_1})({x+a_2})\dots({x+a_{k}})} = \frac{1}{x} - \frac{{a_1}{a_2}\dots{a_{n}}}{x(x+a_1)(x+a_2)\dots(x+a_{n})}
\end{equation*}
  
\section{Preliminaries}

By modifying Ap\'{e}ry's formula we are able to obtain closed form evaluations for a number of series involving Fibonacci and Lucas numbers. \\

But first let us establish the following notation for Binet forms of the Fibonacci and Lucas numbers:

\begin{equation*}
F_n = \frac{\alpha^n - \beta^n}{\sqrt{5}}, \quad L_n = \alpha^n + \beta^n
\end{equation*}

where $$\alpha = \frac{1+\sqrt{5}}{2}, \quad \beta = \frac{1-\sqrt{5}}{2}$$

and $F_n$, $L_n$  are the $n$-th Fibonacci and Lucas number respectively.\\

In addition, we will need the following lemmas.

\begin{lemma} For $n > 0$,

\begin{equation}
1+\alpha^{2n} = \begin{cases}
\sqrt{5}F_n\alpha^n& \text{for odd $n$}\\
L_n\alpha^n& \text{for even $n$}
\end{cases} \label{lem:1}
\end{equation}

\end{lemma}

 This lemma can be found as Exercise 49 on page 97 of Koshy \cite{koshy01}, and may be verified using Binet's formula for the Fibonacci and Lucas numbers.
 
 \begin{lemma} For $m > 0$ and $n \geq 1$,
 
 \begin{equation}
 F_{2^nm}=F_mL_mL_{2m}L_{4m} \cdots L_{2^{n-1}m} \label{lem:2}
 \end{equation}
 
 \begin{proof}
 
 Using the well known identity: $F_{2nm}=F_{nm}L_{nm}$ we find that:
 
 $$F_{2^nm}=F_{2^{n-1}m}L_{2^{n-1}m}=F_{2^{n-2}m}L_{2^{n-2}m}L_{2^{n-1}m}$$
 
 Repeating this process it becomes clear that: 
 
 $$F_{2^nm}=F_mL_mL_{2m}L_{4m} \cdots L_{2^{n-1}m}$$
 
 \end{proof}
 \end{lemma}
 
 \begin{lemma} For any positive integers $q$ and $m$,
 
 \begin{equation}
 F_{(2m+1)q} = F_q\left[(-1)^{mq} + \sum_{k=0}^{m-1}{(-1)^{kq}L_{2(m-k)q}}\right] \label{lem:3.1}
 \end{equation}
 
  \begin{equation}
  L_{(2m+1)q} = L_q\left[(-1)^{m(q+1)} + \sum_{k=0}^{m-1}{(-1)^{k(q+1)}L_{2(m-k)q}}\right] \label{lem:3.2}
  \end{equation}
  
  \begin{proof}
  Only proof of \eqref{lem:3.1} will be given, \eqref{lem:3.2} follows an almost identical procedure.\\ 
  
  Case 1: $q$ odd\\
  
  If $q$ is odd then $(-1)^{jq} = (-1)^{j}$.  In addition, from Binet's formulae for the Fibonacci and Lucas numbers, along with the formula for the finite Geometric Series, we have:\\

    \begin{align*}
     F_{(2m+1)q} &= F_q\left[(-1)^{m} + \sum_{k=0}^{m-1}{(-1)^{k}L_{2(m-k)q}}\right] \\    
      &= \frac{\alpha^{q} - \beta^{q}}{\sqrt{5}}\left[(-1)^{m} +  \sum_{k=0}^{m-1}{(-1)^{k}\left(\alpha^{2(m-k)q} + \beta^{2(m-k)q}\right)}\right] \\   
      &= \frac{\alpha^{q} - \beta^{q}}{\sqrt{5}}\left[(-1)^m +  {\alpha^{2mq}}\sum_{k=0}^{m-1}{(-\alpha^{-2q})^k} + {\beta^{2mq}}\sum_{k=0}^{m-1}{(-\beta^{-2q})^k}\right] \\
      &= \frac{\alpha^{q} - \beta^{q}}{\sqrt{5}}\left[(-1)^m + \frac{\alpha^{2mq}(1-(-1)^{m}\alpha^{-2mq})}{1+\alpha^{-2q}} + \frac{\beta^{2mq}(1-(-1)^{m}\beta^{-2mq})}{1+\beta^{-2q}}\right] \\
      &= \frac{-\beta^q(\alpha^{2q} + 1)}{\sqrt{5}}\left[\frac{(-1)^m(1+\alpha^{2q})}{1+\alpha^{2q}} + \frac{\alpha^{2mq+2q}-(-1)^{m}\alpha^{2q}}{1+\alpha^{2q}} + \frac{\beta^{2mq}-(-1)^{m}}{1+\alpha^{2q}}\right] \\
      &= \frac{-\beta^q}{\sqrt{5}}\left(\alpha^{2mq+2q} + \beta^{2mq}\right) = \frac{\alpha^{2mq+q} - \beta^{2mq+q}}{\sqrt{5}} = F_{(2m+1)q}
    \end{align*}
    
    Case 2: $q$ even\\
    
    If $q$ is even then $(-1)^{jq} = 1$.\\
    
    Therefore, we find that:
    
     \begin{align*}
       F_{(2m+1)q} &= F_q\left[1 + \sum_{k=0}^{m-1}{L_{2(m-k)q}}\right] \\    
        &= \frac{\alpha^{q} - \beta^{q}}{\sqrt{5}}\left[1 +  \sum_{k=0}^{m-1}{\left(\alpha^{2(m-k)q} + \beta^{2(m-k)q}\right)}\right] \\   
        &= \frac{\alpha^{q} - \beta^{q}}{\sqrt{5}}\left[1 +  {\alpha^{2mq}}\sum_{k=0}^{m-1}{(\alpha^{-2q})^k} + {\beta^{2mq}}\sum_{k=0}^{m-1}{(\beta^{-2q})^k}\right] \\
        &= \frac{\alpha^{q} - \beta^{q}}{\sqrt{5}}\left[1 + \frac{\alpha^{2mq}(1-\alpha^{-2mq})}{1-\alpha^{-2q}} + \frac{\beta^{2mq}(1-\beta^{-2mq})}{1-\beta^{-2q}}\right] \\
        &= \frac{-\beta^q(1-\alpha^{2q})}{\sqrt{5}}\left[\frac{(1-\alpha^{2q})}{1-\alpha^{2q}} + \frac{\alpha^{2q}-\alpha^{2mq+2q}}{1-\alpha^{2q}} + \frac{\beta^{2mq}-1}{1-\alpha^{2q}}\right] \\
        &= \frac{-\beta^q}{\sqrt{5}}\left(\beta^{2mq} - \alpha^{2mq+2q}\right) = \frac{\alpha^{2mq+q} - \beta^{2mq+q}}{\sqrt{5}} = F_{(2m+1)q}
      \end{align*}

  \end{proof}
  
  \begin{remark}
  As an interesting side note, by Lemma 3 we have:
  
  \begin{align*}
   F_{(2l+1)(2m+1)} &= F_{2m+1}\left[(-1)^l + \sum_{k=0}^{l-1}{(-1)^{k}L_{2(l-k)(2m+1)}}\right]\\
     &= F_{2l+1}\left[(-1)^m + \sum_{k=0}^{m-1}{(-1)^{k}L_{2(m-k)(2l+1)}}\right]
  \end{align*}
  
  So that for positive integers $l$ and $m$,
  
  \begin{equation*}
  \frac{F_{2l+1}}{F_{2m+1}} = \frac{(-1)^l + \sum_{k=0}^{l-1}{(-1)^{k}L_{2(l-k)(2m+1)}}}{(-1)^m + \sum_{k=0}^{m-1}{(-1)^{k}L_{2(m-k)(2l+1)}}}
  \end{equation*}
  
  and,
  
  \begin{equation*}
  \frac{L_{2l+1}}{L_{2m+1}} = \frac{1 + \sum_{k=0}^{l-1}{L_{2(l-k)(2m+1)}}}{1 + \sum_{k=0}^{m-1}{L_{2(m-k)(2l+1)}}}
  \end{equation*}
  
  \end{remark}
  
 \end{lemma}
 
 The following lemma is a generalization of two identities of Usiskin \cite{usiskin74a}, \cite{usiskin74b}; who proposed them as problems in the \textit{Fibonacci Quarterly} in 1974.  In \cite{filipponi96} Filipponi also achieved similar results to those presented here.  In fact his \textit{Proposition 2} (presented here as Lemma 7) is equivalent to \eqref{lem:4.1}.  Although, to the best of the author's knowledge equation \eqref{lem:4.2} is a novel result.
 
 \begin{lemma} For odd $p$, $F_{p^n}$ and $L_{p^n}$ may be written as:

\begin{equation}
F_{p^n} = F_{(2m+1)^n} = \prod_{j=0}^{n-1}\left[{(-1)^m + \sum_{k=0}^{m-1}{(-1)^k{L_{2(m-k)(2m+1)^j}}}}\right] \label{lem:4.1}
\end{equation}

\begin{equation}
L_{p^n} = L_{(2m+1)^n} = \prod_{j=0}^{n-1}\left[{1 + \sum_{k=0}^{m-1}{{L_{2(m-k)(2m+1)^j}}}}\right] \label{lem:4.2}
\end{equation}

\begin{proof}
Using \eqref{lem:3.1} and \eqref{lem:3.2} from Lemma 3 we are able to prove Lemma 4.  We proceed by induction on $n$.\\

Base Case: $n=1$

\begin{align*}
F_1\prod_{j=0}^{0}\left[{(-1)^m + \sum_{k=0}^{m-1}{(-1)^k{L_{2(m-k)(2m+1)^j}}}}\right] &= F_1\left[{(-1)^m + \sum_{k=0}^{m-1}{(-1)^k{L_{2(m-k)}}}}\right] \\
  &= F_{2m+1} 
\end{align*}

where we have used \eqref{lem:3.1} with $q=1$.\\

Induction Step:\\

Assume that,  

\begin{equation*}
F_{p^n} = F_{(2m+1)^n} = \prod_{j=0}^{n-1}\left[{(-1)^m + \sum_{k=0}^{m-1}{(-1)^k{L_{2(m-k)(2m+1)^j}}}}\right]
\end{equation*}

We now prove that \eqref{lem:4.1} holds for $n+1$.\\

That is,

\begin{align*}
\prod_{j=0}^{n}\left[{(-1)^m + \sum_{k=0}^{m-1}{(-1)^k{L_{2(m-k)(2m+1)^j}}}}\right] &= F_{(2m+1)^n}\left[{(-1)^m + \sum_{k=0}^{m-1}{(-1)^k{L_{2(m-k)(2m+1)^n}}}}\right] \\
&=  F_{(2m+1)^{n+1}}
\end{align*}

where we have used \eqref{lem:3.1} with $q = (2m+1)^n$.\\

This concludes the proof of equation \eqref{lem:4.1}, the proof of equation \eqref{lem:4.2} proceeds in a nearly identical fashion.
\end{proof}

\end{lemma}

The following four lemmas are due to Filipponi, and can found in \cite{filipponi96}; as a result, they are presented without proof. Although it is not too difficult to see that the proofs of Lemmas 5-8 follow closely to that of Lemma 4, in that heavy use is made of both Binet's formulae for the Fibonacci and Lucas numbers as well as the formula for the finite geometric series.  

\begin{lemma}
For even $p$; and $n, m \geq1$,

\begin{equation}
F_{mp^n} = F_{mp}\prod_{j=1}^{n-1}{\left[\sum_{k=1}^{p/2}{L_{(2k-1)mp^j}}\right]}
\end{equation}

\end{lemma}

\begin{lemma}
For odd $p$, even $m$ and $n\geq1$,

\begin{equation}
F_{mp^n} = F_{mp}\prod_{j=1}^{n-1}{\left[1+\sum_{k=1}^{(p-1)/2}{L_{2mp^j}}\right]}
\end{equation}

\end{lemma}

\begin{lemma}
For odd $p$, and $n\geq1$,

\begin{equation}
F_{p^n} = (-1)^{(n-1)(p-1)/2}F_{p}\prod_{j=1}^{n-1}{\left[1+\sum_{k=1}^{(p-1)/2}{(-1)^k{L_{2kp^j}}}\right]}
\end{equation}

\end{lemma}

\begin{lemma}
For even $p$ and $n \geq 2$,

\begin{equation}
L_{p^n} = 2 + (L_{p^2}-2)\prod_{j=1}^{n-2}{\left[\sum_{k=1}^{p/2}{L_{(2k-1)p^{j+1}/2}}\right]^2}
\end{equation}

\end{lemma}
 
 \section{Results}

In 1974 D.A. Millin \cite{millin74} proposed a problem in \textit{The Fibonacci Quarterly} regarding the series of reciprocal Fibonacci numbers of the form $F_{2^n}$ (coincidentally in the same issue that Usiskin's problems appeared). Millin's proposal stimulated significant interest in the series, which led to the publication of numerous proofs and generalizations (see \cite{hoggatt76.3} and \cite{melham95} for more information and references).\\

For the first theorem, a new proof of Millin's Series is given using an Ap\'{e}ry-like formula.

\begin{theorem}

\begin{equation}
\sum_{n=0}^{\infty}{\frac{1}{F_{2^n}}}=\frac{7-\sqrt{5}}{2}
\end{equation}

\begin{proof}

Let,

\begin{equation*}
B_n = \frac{1}{(1+\alpha^{4})(1+\alpha^{8})\cdots (1+\alpha^{2^{n}})}, \quad n > 1
\end{equation*}

By subtracting $B_{n+1}$ from $B_n$ we observe that:

\begin{equation}
B_n - B_{n+1} = \frac{\alpha^{2^{n+1}}}{(1+\alpha^{4})(1+\alpha^{8})\cdots (1+\alpha^{2^{n}})(1+\alpha^{2^{n+1}})} \label{Bseq}
\end{equation}

Setting $n = 2^n$ in Lemma 1 we obtain:

$$1+\alpha^{2^{n+1}} = L_{2^n}\alpha^{2^n}$$

Note that we are not able to begin $B_n$ at $(1+\alpha^2)^{-1}$ since the above altered form of Lemma 1 does not hold for $n = 0$.\\

By Lemma 1, we now find that \eqref{Bseq} becomes:

\begin{equation*}
 B_n - B_{n+1} = \frac{\alpha^{2^{n+1}}}{L_{2}\alpha^{2} L_{4}\alpha^{4}\cdots L_{2^n}\alpha^{2^{n}}} = \frac{\alpha^{2^{n+1}}}{L_{2}L_{4}\cdots L_{2^n}\alpha^{2^{n+1}-2}}
\end{equation*}
 
 Then by cancellation and Lemma 2 we obtain:
 
\begin{equation*}
B_n - B_{n+1} = \frac{\alpha^2}{F_{2^{n+1}}}
\end{equation*}

Using this information, we are now able to evaluate Millin's series telescopically.\\

In other words,

\begin{align*}
\sum_{n=2}^{\infty}\frac{\alpha^2}{F_{2^{n+1}}} &= \sum_{n=2}^{\infty}B_n - B_{n+1} \\
  &= B_2 - \lim_{n \to \infty} B_{n+1} \\
  &= \frac{1}{1+\alpha^4} - 0
\end{align*}

Therefore,

\begin{equation}
\sum_{n=2}^{\infty}\frac{1}{F_{2^{n+1}}} = \sum_{n=3}^{\infty}\frac{1}{F_{2^{n}}} = \frac{1}{\alpha^2 + \alpha^6}
\end{equation}

Since Millin's series begins the sum at $n=0$, we must add the missing terms.\\

In the end we find that: 

\begin{equation}
\sum_{n=0}^{\infty}\frac{1}{F_{2^{n}}} = \frac{1}{\alpha^2+\alpha^6} + 1 + 1 + 1/3 = \frac{7-\sqrt{5}}{2}
\end{equation}

\end{proof}

\end{theorem}

\begin{theorem} For any positive integer $m$, and $a > 0$,

\begin{equation}
\sum_{n=0}^{\infty}{\frac{L^a_{2^{n+1}m}-1}{F^a_{2^{n+2}m}}} = \frac{1}{F^a_mL^a_m} \label{thrm:2}
\end{equation}

\begin{proof}

Let,

\begin{equation*}
B_0 = \frac{1}{F^a_mL^a_m}, \quad B_n = \frac{1}{F^a_mL^a_mL^a_{2m}L^a_{4m}\cdots{L^a_{{2^n}m}}} \quad n > 0
\end{equation*} 

Then subtracting $B_{n+1}$ from $B_n$ and employing Lemma 2 gives,

\begin{equation*}
B_{n} - B_{n+1} = \frac{L^a_{2^{n+1}m}-1}{F^a_{2^{n+2}m}}
\end{equation*}

Once again, taking advantage of the telescoping, we acquire:

\begin{equation*}
\sum_{n=0}^{\infty}{\frac{L^a_{2^{n+1}m}-1}{F^a_{2^{n+2}m}}} = B_0 - \lim_{n \to \infty}B_{n+1} =\frac{1}{F^a_mL^a_m}
\end{equation*}

\end{proof}

\begin{remark}
Specifically, setting $a=m=1$, \eqref{thrm:2} becomes:

\begin{equation}
\sum_{n=0}^{\infty}{\frac{L_{2^{n+1}}-1}{F_{2^{n+2}}}} = 1 
\end{equation}

Although since we know the closed form for Millin's Series, we may split up the summand to obtain,

\begin{equation*}
\sum_{n=0}^{\infty}{\frac{L_{2^{n+1}}}{F_{2^{n+2}}}} = 1 + \sum_{n=0}^{\infty}{\frac{1}{F_{2^{n+2}}}} = \frac{5-\sqrt{5}}{2}
\end{equation*}

\end{remark}

\end{theorem}

\begin{theorem} For any positive integer $m$,

\begin{equation}
\sum_{n=0}^{\infty}{\frac{{(-1)^m - 1 + \sum_{k=0}^{m-1}{(-1)^k{L_{2(m-k)(2m+1)^n}}}}}{F_{(2m+1)^{n+1}}}} = 1 \label{thrm:3.1}
\end{equation}

\begin{proof}

Let,

\begin{equation*}
B_0 = 1, \quad B_n = \frac{1}{\prod_{j=0}^{n-1}\left[{(-1)^m + \sum_{k=0}^{m-1}{(-1)^k{L_{2(m-k)(2m+1)^j}}}}\right]}, \quad n > 0
\end{equation*}

Then,

\begin{align*}
B_n - B_{n+1} &= \frac{\left[{(-1)^m + \sum_{k=0}^{m-1}{(-1)^k{L_{2(m-k)(2m+1)^n}}}}\right] - 1}{\prod_{j=0}^{n}\left[{(-1)^m + \sum_{k=0}^{m-1}{(-1)^k{L_{2(m-k)(2m+1)^j}}}}\right]} \\
\text{which by Lemma 4 gives,} \\
  &= \frac{{(-1)^m - 1 + \sum_{k=0}^{m-1}{(-1)^k{L_{2(m-k)(2m+1)^n}}}}}{F_{(2m+1)^{n+1}}}
\end{align*}

Therefore, 

\begin{equation*}
\sum_{n=0}^{\infty}{\frac{{(-1)^m - 1 + \sum_{k=0}^{m-1}{(-1)^k{L_{2(m-k)(2m+1)^n}}}}}{F_{(2m+1)^{n+1}}}} = B_0 - \lim_{n \to \infty} B_{n+1} = 1
\end{equation*}

as desired.\\

\end{proof}

\begin{remark}
As a specific example, let $m=1$. \\

We then have,

\begin{equation}
\sum_{n=0}^{\infty}{\frac{L_{2\cdot3^{n}}-2}{F_{3^{n+1}}}} = \sum_{n=0}^{\infty}{\frac{L^2_{3^{n}}}{F_{3^{n+1}}}} = 1 \label{thrm:3.2}
\end{equation}

where we have used the identity $L^2_n = L_{2n} + 2(-1)^n$ \cite[page 97]{koshy01}.\\

In addition, by splitting up the summand of sum on the left side of \eqref{thrm:3.2}, we achieve a series representation for the series of reciprocal Fibonacci numbers $F_{3^n}$.\\

Namely,

\begin{equation}
\sum_{n=0}^{\infty}{\frac{1}{F_{3^{n}}}} = \frac{1}{2}\sum_{n=0}^{\infty}{\frac{L_{2\cdot3^{n}}}{F_{3^{n+1}}}} + \frac{1}{2}
\end{equation}

More generally, for all odd $m$ \eqref{thrm:3.1} becomes:

\begin{equation*}
\sum_{n=0}^{\infty}{\frac{{\left[\sum_{k=0}^{m-1}{(-1)^k{L_{2(m-k)(2m+1)^n}}}\right]-2}}{F_{(2m+1)^{n+1}}}} = 1 
\end{equation*}

Therefore for all Fibonacci numbers of the form $F_{4m+3}$,

\begin{equation}
\sum_{n=0}^{\infty}{\frac{1}{F_{(4m+3)^{n}}}} = \frac{1}{2}\sum_{n=0}^{\infty}{\frac{{\sum_{k=0}^{2m}{(-1)^k{L_{2(2m+1-k)(4m+3)^n}}}}}{F_{(4m+3)^{n+1}}}} + \frac{1}{2} 
\end{equation}

For instance,

\begin{equation}
\sum_{n=0}^{\infty}{\frac{1}{F_{7^{n}}}} = \frac{1}{2}\sum_{n=0}^{\infty}{\frac{L_{6\cdot7^{n}} - L_{4\cdot7^{n}} + L_{2\cdot7^{n}}}{F_{7^{n+1}}}} + \frac{1}{2} 
\end{equation}

\end{remark}

\end{theorem}

\begin{theorem} For any positive integer $m$,

\begin{equation}
\sum_{n=0}^{\infty}{\frac{{\sum_{k=0}^{m-1}{{L_{2(m-k)(2m+1)^n}}}}}{L_{(2m+1)^{n+1}}}} = 1 \label{thrm:4.1}
\end{equation}

\begin{proof}
In the interest of brevity, a full proof of Theorem 4 will be omitted, as it follows the same method of proof as Theorem 3, except for the use of \eqref{lem:4.2} of Lemma 4 instead of \eqref{lem:4.1}.

\end{proof}

\begin{remark}
To illustrate \eqref{thrm:4.1} with a particular case, let $m=1$. \\ 

We obtain the series:

 \begin{equation}
 \sum_{n=0}^{\infty}{\frac{L_{2\cdot3^{n}}}{L_{3^{n+1}}}} = 1 \label{thrm:4.2}
 \end{equation}

\end{remark}

\end{theorem}

\begin{theorem} For any positive integer $m$,

\begin{equation}
\sum_{n=0}^{\infty}\frac{{F_{2^{n+2}}}\left[(-1)^{m} + \sum_{k=0}^{m-1}{(-1)^{k}L_{(m-k)2^{j+1}}}\right]}{F_{(2m+1)2^{n+2}}} = \frac{1}{F_{2m+1}L_{2m+1}} \label{thrm:5.1}
\end{equation}

\begin{proof}

Setting  $q=2^n$ in Lemma 3 returns:

\begin{equation}
(-1)^{m} + \sum_{k=0}^{m-1}{(-1)^{k}L_{(m-k)2^{n+1}}} = \begin{cases}
F_{2m+1}, \quad \text{if $n=0$}\\
\frac{L_{(2m+1)2^n}}{L_{2^n}}, \quad \text{if $n>0$}
\end{cases} \label{thrm:5.2}
\end{equation}

Then, setting $B_n$ to:

\begin{equation*}
B_n = \frac{1}{\prod_{j=0}^{n}\left[(-1)^{m} + \sum_{k=0}^{m-1}{(-1)^{k}L_{(m-k)2^{j+1}}}\right]}, \quad n \geq 0
\end{equation*}

we have for $B_0$,

\begin{equation*}
B_0 = \frac{1}{\left[(-1)^{m} + \sum_{k=0}^{m-1}{(-1)^{k}L_{(m-k)}}\right]} = \frac{1}{F_{2m+1}}
\end{equation*}

In addition, by subtracting $B_{n+1}$ from $B_n$ produces:

\begin{align*}
B_n - B_{n+1} &= \frac{\left[(-1)^{m} + \sum_{k=0}^{m-1}{(-1)^{k}L_{(m-k)2^{n+2}}}\right]-1}{\prod_{j=0}^{n+1}\left[(-1)^{m} + \sum_{k=0}^{m-1}{(-1)^{k}L_{(m-k)2^{j+1}}}\right]}\\ 
\text{which by \eqref{thrm:5.2} gives:}\\
  &= \frac{\left[(-1)^{m} -1 + \sum_{k=0}^{m-1}{(-1)^{k}L_{(m-k)2^{j+2}}}\right]\prod_{j=1}^{n+1}{L_{2^j}}}{F_{2m+1}\prod_{j=1}^{n+1}{L_{(2m+1)2^j}}}\\
  &= \frac{L_{2m+1}\left[(-1)^{m} -1 + \sum_{k=0}^{m-1}{(-1)^{k}L_{(m-k)2^{j+2}}}\right]F_1\prod_{j=0}^{n+1}{L_{2^j}}}{F_{2m+1}\prod_{j=0}^{n+1}{L_{(2m+1)2^j}}} \\
  &= \frac{L_{2m+1}{F_{2^{n+2}}}\left[(-1)^{m} -1 + \sum_{k=0}^{m-1}{(-1)^{k}L_{(m-k)2^{j+2}}}\right]}{F_{(2m+1)2^{n+2}}}
\end{align*}

Thus,

\begin{equation*}
\sum_{n=0}^{\infty}\frac{L_{2m+1}{F_{2^{n+2}}}\left[(-1)^{m} -1 + \sum_{k=0}^{m-1}{(-1)^{k}L_{2(m-k)2^{n+1}}}\right]}{F_{(2m+1)2^{n+2}}} = B_0 = \frac{1}{F_{2m+1}} 
\end{equation*}

and,

\begin{equation*}
\sum_{n=0}^{\infty}\frac{{F_{2^{n+2}}}\left[(-1)^{m} -1 + \sum_{k=0}^{m-1}{(-1)^{k}L_{2(m-k)2^{n+1}}}\right]}{F_{(2m+1)2^{n+2}}} =  \frac{1}{F_{2m+1}L_{2m+1}} 
\end{equation*}

\end{proof}

\begin{remark}

As an example of a specific case, let $m=1$. \\

This gives,

 \begin{equation}
 \sum_{n=0}^{\infty}{\frac{F_{2^{n+3}}-2F_{2^{n+2}}}{F_{3\cdot2^{n+2}}}} = \frac{1}{8} 
 \end{equation}

\end{remark}

\end{theorem}

\begin{theorem} For even $p$ and $m \geq 1$,

\begin{equation}
\sum_{n=1}^{\infty}\frac{\left[ \sum_{k=1}^{p/2}{L_{(2k-1)mp^{n}}}\right] -1}{F_{mp^{n+1}}} = \frac{1}{F_{mp}} \label{thrm:6.1}
\end{equation}

\begin{proof}
As before, let

\begin{equation*}
B_1 = 1, \quad B_n = \frac{1}{\prod_{j=1}^{n-1}{\left[\sum_{k=1}^{p/2}{L_{(2k-1)mp^j}}\right]}}, \quad n > 1
\end{equation*}

As a result, subtracting $B_{n+1}$ from $B_n$ we obtain,

\begin{equation*}
B_{n}-B_{n+1} = \frac{\left[\sum_{k=1}^{p/2}{L_{(2k-1)mp^n}}\right]-1}{\prod_{j=1}^{n}{\left[\sum_{k=1}^{p/2}{L_{(2k-1)mp^j}}\right]}}
\end{equation*}

which by Lemma 5 gives:

\begin{equation}
B_{n}-B_{n+1} = \frac{F_{mp}\left[\left(\sum_{k=1}^{p/2}{L_{(2k-1)mp^n}}\right)-1\right]}{F_{mp^{n+1}}}
\end{equation}

Therefore, we have:

\begin{equation}
\sum_{n=1}^{\infty}\frac{F_{mp}\left[\left(\sum_{k=1}^{p/2}{L_{(2k-1)mp^n}}\right)-1\right]}{F_{mp^{n+1}}} = B_1 = 1
\end{equation}

which becomes \eqref{thrm:6.1} after both sides are divided by $F_{mp}$.

\end{proof}

\begin{remark}
Just as Theorem 3 provides us with a series representation for reciprocal Fibonacci numbers of the form $F_{4m+3}$, we are now able to give a series representation for the reciprocals of even indexed Fibonacci numbers.\\

Namely, even $p$ and $m\geq1$,

\begin{equation}
\sum_{n=1}^{\infty}{\frac{1}{F_{mp^{n+1}}}} = \left[\sum_{n=1}^{\infty}\frac{\sum_{k=1}^{p/2}{L_{(2k-1)mp^n}}}{F_{mp^{n+1}}}\right] - \frac{1}{F_{mp}}
\end{equation}

\end{remark}

\end{theorem}

\begin{theorem} For odd $p$ and even $m$,

\begin{equation}
\sum_{n=1}^{\infty}\frac{\sum_{k=1}^{(p-1)/2}{L_{2kmp^n}}}{F_{mp^{n+1}}} = \frac{1}{F_{mp}} \label{thrm:7.1}
\end{equation}

\begin{proof}
A proof of Theorem 7 will be omitted since it follows the same steps as Theorem 6, except for the use of Lemma 6 instead of Lemma 5.
\end{proof}

\begin{remark}
To demonstrate \eqref{thrm:7.1}, let $p = 5$ and $m = 4$.  We then have,

\begin{equation}
\sum_{n=1}^{\infty}\frac{L_{8\cdot5^n} + L_{16\cdot5^n}}{F_{4\cdot5^{n+1}}} = \frac{1}{6765} 
\end{equation}

\end{remark}

\end{theorem}

\begin{theorem}
For odd $p$,

\begin{equation}
\sum_{n=1}^{\infty}{\frac{(-1)^{n(p-1)/2}\sum_{k=1}^{(p-1)/2}{(-1)^k{L_{2kp^n}}}}{F_{p^{n+1}}}} = \frac{1}{F_p} \label{thrm:8.1}
\end{equation}

\begin{proof}
Let,

\begin{equation*}
B_1 = 1, \quad B_n = \frac{1}{\prod_{j=1}^{n-1}{\left[1+\sum_{k=1}^{(p-1)/2}{(-1)^k{L_{2kp^j}}}\right]}}, \quad n > 1
\end{equation*}

Therefore,

\begin{equation*}
B_n - B_{n+1} = \frac{\sum_{k=1}^{(p-1)/2}{(-1)^k{L_{2kp^n}}}}{\prod_{j=1}^{n}{\left[1+\sum_{k=1}^{(p-1)/2}{(-1)^k{L_{2kp^j}}}\right]}}
\end{equation*}

Then utilizing Lemma 7, we have:

\begin{equation}
B_n - B_{n+1} = \frac{(-1)^{n(p-1)/2}F_p\left[\sum_{k=1}^{(p-1)/2}{(-1)^k{L_{2kp^n}}}\right]}{F_{p^{n+1}}}
\end{equation}

which after summation produces \eqref{thrm:8.1}.

\end{proof}

\begin{remark}
As an example of \eqref{thrm:8.1}, for $p=3$ we have:

\begin{equation}
\sum_{n=1}^{\infty}{\frac{(-1)^{n+1}{L_{2\cdot3^n}}}{F_{3^{n+1}}}} = \frac{1}{2}
\end{equation}

\end{remark}

\end{theorem}

\begin{theorem} For even $p$,

\begin{equation}
\sum_{n=2}^{\infty}{\frac{\left[\sum_{k=1}^{p/2}{L_{(2k-1)p^n/2}}\right]^2-1}{L_{p^{n+1}}-2}} = \frac{1}{L_{p^2} - 2} \label{thrm:9.1}
\end{equation}

\begin{proof}
Let,

\begin{equation*}
B_2 = 1, \quad B_n = \frac{1}{\prod_{j=1}^{n-2}{\left[\sum_{k=1}^{p/2}{L_{(2k-1)p^{j+1}/2}}\right]^2}}, \quad n > 2
\end{equation*}

we then have:

\begin{equation*}
B_n - B_{n+1} = \frac{\left[\sum_{k=1}^{p/2}{L_{(2k-1)p^{n}/2}}\right]^2-1}{\prod_{j=1}^{n-1}{\left[\sum_{k=1}^{p/2}{L_{(2k-1)p^{j+1}/2}}\right]^2}}
\end{equation*}

which by Lemma 8 gives,

\begin{equation*}
B_n - B_{n+1} = \frac{(L_{p^2}-2)\left[\left(\sum_{k=1}^{p/2}{L_{(2k-1)p^{n}/2}}\right)^2-1\right]}{L_{p^{n+1}-2}}
\end{equation*}

Once again, by carrying out the summation, we find that:

\begin{equation*}
\sum_{n=2}^{\infty}{\frac{\left[\sum_{k=1}^{p/2}{L_{(2k-1)p^n/2}}\right]^2-1}{L_{p^{n+1}}-2}} = \frac{B_2}{L_{p^2} - 2} = \frac{1}{L_{p^2} - 2}
\end{equation*}

\end{proof}

\begin{remark}
To demonstrate \eqref{thrm:9.1}, let $p = 2$. \\  

We then have:

\begin{equation}
\sum_{n=2}^{\infty}{\frac{L^2_{2^{n-1}}-1}{L_{2^{n+1}}-2}} = \frac{1}{5}
\end{equation}

\end{remark}

\end{theorem}

\begin{closeremark}
The aim of this paper was to demonstrate how Ap\'{e}ry-like formulae may be employed to derive certain infinite series involving Fibonacci and Lucas numbers.  While this method is not well suited for all series of this type, it provides us with another useful tool to be used in the evaluation of series involving Fibonacci and Lucas numbers.

In addition, over the course of this paper we have restricted ourselves to examining series involving Fibonacci and Lucas numbers of the form $F_{p^k}$ and $L_{p^k}$.  This due to the fact that terms of this form usually yield "clean" looking summands.  That being said, these are not the only series that may be evaluated using Ap\'{e}ry-like formulae.  Series of other types may also considered, although they may not always return such pretty results.

Lastly, the techniques employed in this paper should yield similar dividends when applied to more general recursively-defined sequences, such as the Generalized Fibonacci Sequences $U_n$ and $V_n$.
\end{closeremark}

\begin{acknolwedgment}
The research for, as well as the writing of this paper, was carried out over the Spring and Summer of 2015 as part of an undergraduate research project. The author would like to express his sincerest gratitude to Dr. Al Lehnen; whose suggestions, corrections, and contributions proved to be invaluable, and greatly enhanced the quality of this paper.
\end{acknolwedgment}

MSC2010: 11B39\\

Madison College, Madison, Wisconsin. [Student]\\

\textit{E-mail: sanfordchance@gmail.com}
\end{document}